\documentclass[12pt]{amsart}

\usepackage{verbatim}

\usepackage{tikz-cd}

\usepackage{mathabx}

\usepackage[backref,bookmarks=false]{hyperref}

\usepackage{color}

\newcommand{\lab}{\label}

\numberwithin{equation}{section}

\newcommand{\ben}{\begin{enumerate}}
\newcommand{\een}{\end{enumerate}}

\newcommand{\bea}{\begin{eqnarray}}
\newcommand{\ba}{\begin{array}}
\newcommand{\bean}{\begin{eqnarray*}}
\newcommand{\ea}{\end{array}}
\newcommand{\eea}{\end{eqnarray}}
\newcommand{\eean}{\end{eqnarray*}}
\newcommand{\beq}{\begin{equation}}
\newcommand{\eeq}{\end{equation}}
\newcommand{\bthm}{\begin{thm}}
\newcommand{\ethm}{\end{thm}}
\newcommand{\blem}{\begin{lem}}
\newcommand{\elem}{\end{lem}}
\newcommand{\bprop}{\begin{prop}}
\newcommand{\eprop}{\end{prop}}
\newcommand{\bcor}{\begin{cor}}
\newcommand{\ecor}{\end{cor}}
\newcommand{\bdfn}{\begin{dfn}}
\newcommand{\edfn}{\end{dfn}}
\newcommand{\brem}{\begin{rem}}
\newcommand{\erem}{\end{rem}}
\newcommand{\bpf}{\begin{proof}}
\newcommand{\epf}{\end{proof}}
\newcommand{\bfact}{\begin{fact}}
\newcommand{\efact}{\end{fact}}
\newcommand{\bobs}{\begin{obs}}
\newcommand{\eobs}{\end{obs}}
\newcommand{\bexam}{\begin{exam}}
\newcommand{\eexam}{\end{exam}}
\newcommand{\bclaim}{\begin{claim}}
\newcommand{\eclaim}{\end{claim}}

\newtheorem{thm}{Theorem}[section]
\newtheorem{prop}[thm]{Proposition}
\newtheorem{lem}[thm]{Lemma}

\newtheorem{cor}[thm]{Corollary}
\newtheorem{dfn}[thm]{Definition}
\newtheorem{rem}[thm]{Remark}
\newtheorem{fact}[thm]{Fact}
\newtheorem{claim}[thm]{Claim}
\newtheorem{obs}[thm]{Observation}
\newtheorem{exam}[thm]{Example}

\newtheorem*{condition'}{Condition 2'}

 \newtheoremstyle{claimstyle}%
   {}
   {}
   {\normalfont}
   {}
   {\itshape}
   {.}
   { }
   {\thmnote{#3}}

\theoremstyle{claimstyle}

\alph{enumii} \roman{enumiii}

             \def\cB{\mathcal B}       \def\cC{\mathcal C}
\def\cH{\mathcal H}                    
                   
                       \newcommand{\J}{\mathcal{J}}
                           
           \def\cK{\mathcal K}

                      \def\R{{\mathbb R}}
\def\C{{\mathbb C}}

                             \def\d{\delta}
               \def\e{\varepsilon}           
\def\g{\gamma}                            
                         \def\Om{\Omega}
               \def\sg{\sigma}
                          
\def\ka{\kappa}

\newcommand{\ph}{\varphi}
\newcommand{\al}{\alpha}
\newcommand{\ga}{\gamma}

\def\1{1\!\!1}

\def\and{\text{ and }}

     \def\HD{\text{{\rm HD}}}   
         
\def\HypDim{\text{{\rm HypDim}}}

\def\Int{\text{{\rm Int}}}
         \def\P{\text{{\rm P}}}

\def\({\bigl(}                \def\){\bigr)}
\def\lt{\left}                \def\rt{\right}

                        \def\^{\tilde}

            \def\sms{\setminus}
\def\sbt{\subset}

\def\sp{\medskip}                     
\def\ov{\overline}

\def\arg{\text{arg}}

\def\D{{\mathbb D}}

\def\${$ \displaystyle }

\def\disf{{\bf f}}
\def\disF{{\bf E}}
\def\disE{{\bf E}}

\newcommand{\pf}{{\mathcal{L}}}

\newcommand{\jul}{\mathcal J}

\def\hypdim{\HypDim(f)}
\def\hypd{\HypDim}

\def\Tract{\Omega}

\def\\HypDim{\HypDim}

\begin{document}

\title[On hyperbolic dimension gap for entire functions]{On hyperbolic dimension gap for entire functions}


\author[\sc Volker MAYER]{\sc Volker Mayer}
\address{Volker Mayer, Universit\'e de Lille, D\'epartement de
  Math\'ematiques, UMR 8524 du CNRS, 59655 Villeneuve d'Ascq Cedex,
  France} \email{volker.mayer@univ-lille.fr \newline
  \hspace*{0.42cm} \it Web: \rm math.univ-lille1.fr/$\sim$mayer}

\author{Mariusz Urba\'nski}
\address{Mariusz Urba\'nski, Department of Mathematics, University of North Texas, Denton, TX 76203-1430, USA}
\email{urbanski@unt.edu \newline \hspace*{0.42cm} \it Web: \rm www.math.unt.edu/$\sim$urbanski}

\date{\today} 
\thanks{The research of Volker Mayer was supported in part by the ANR-DFG project QuaSiDy ANR-21-CE40-0016. The research of Mariusz Urba\'nski was supported in part by the Simons grants 581668 and 900989. Mariusz Urba\'nski also thanks the SRMI in Sydney for warm hospitality and support during the work on this research}


\begin{abstract} Polynomials and entire functions whose hyperbolic dimension is strictly smaller
than the Hausdorff dimension of their Julia set are known to exist but in all these examples the latter dimension
is maximal, i.e. equal to two. In this paper we show that there exist hyperbolic entire functions $f$
having Hausdorff dimension of the Julia set $\HD (\J _f)<2$ and  hyperbolic dimension
$\HypDim(f)<\HD(\J_f)$.
 \end{abstract}

\maketitle



\section{Introduction} \label{s1}

In this paper we consider some relations
between the Hausdorff dimension $\HD (\jul _f )$  and  the hyperbolic dimension $\hypdim $ of the Julia set $\jul_f$ of an entire function 
$f:\C\to\C$,
where $\C$, as usually, denotes the complex plane. More precisely, we show the following.

\bthm\label{thm main}
There exist hyperbolic entire functions $f$ in the Eremenko--Lyubich class $\cB$ such that
$$
\hypdim <\HD (\jul _f )<2.
$$
\ethm

The concept $\hypdim $ of hyperbolic dimension has been introduced by Shishikura  in \cite{Shishi1998}. Given an entire function 
$f:\C\to\C$ it is defined to be the supremum of Hausdorff dimensions of all hyperbolic sets of $f$. We recall that a set $X\sbt\C$ is hyperbolic if it is compact, forward--invariant under $f$ and if there exist $c>0$ and $\kappa >1$ such that
$$
|(f^n)'(z)|\geq c \kappa^n \quad \text{for every $z\in X$ and all $n\geq1$.}
$$
It is immediate from this definition that $X\sbt \jul_f$.
For hyperbolic polynomials the whole Julia set is a hyperbolic set, whence there is no difference between 
the hyperbolic dimension  and  the Hausdorff dimension of the Julia set.
In general, since hyperbolic sets of $f$ are subsets of the Julia set of $f$, we have that
\beq\label{eq 1}\hypdim \leq \HD (\jul _f ).\eeq
Examples of entire functions with strict inequality are known (\cite{St99-1}, \cite{UZ03}). Quite recently Avila-Lyubich \cite{AvLyuI, AvLyuII} showed that there
exist Feigenbaum polynomials having this property. But in all known examples with strict inequality in \eqref{eq 1} the Hausdorff dimension of the Julia set is maximal, i.e. equal to two
and Avila-Lyubich  mention that for arbitrary polynomials $f$ with $\HD (\jul _f )<2$ one should have equality.
Here we show that this is not the case for entire functions even inside the Eremenko-Lyubich class $\cB$ consisting of all entire functions having a bounded  set of finite singularities.


\smallskip

In order to prove Theorem~\ref{thm main} we first need good candidates of entire functions whose Julia sets have Hausdorff dimension 
less than two. 
The first such examples where provided by Gwyneth Stallard during 1990's. The interested reader can find an overview in her survey in \cite{RiSt-Bakerbook}. 
These examples are entire functions having one single logarithmic tract over infinity; see Section 2.1 for the definition of the singularities of entire functions, in particular of logarithmic tract. As nowadays it is well known, the geometry of such a tract
or the growth of the function in the tract influences the size of the Julia set. 
Particularly interesting for the present work is her
family of intermediate growth in \cite{Stallard-III}. The growth does depend on a parameter $p>0$ and these functions are defined by the formula
\beq\label{eq 2}
E(z):= \frac{1}{2i\pi}\int_L \frac{\exp\big( e^{(\log \xi)^{1+p})}\big)}{\xi-z} d\xi\, ,
\eeq
where $L$ is the boundary of the region
\beq\label{eq 3}
G=\left \{x+iy\in\C: |y| < \frac{\pi x }{ (1+p)(\log x)^p} \; , \; x > 3\right\},
\eeq
oriented in the clockwise direction, for $z\in \C\setminus  \ov G$ and by analytic continuation for $z\in \ov G$.
Appropriate details of such analytic extension are given in Section 2.2..
The reader should have in mind that this function is close to 
\beq\label{eq 4}
f(z)= \exp\big( e^{(\log z)^{1+p}}\big)  \quad \text{ for $z\in  G$}
\eeq
and is bounded elsewhere. Here $(\log z)^{1+p}$ is defined so that it gives real values for real $z>e$.

Consider then the family $\(\disF_l:\C\to\C\)_{l\in\C}$ defined by the formula
$$
\disF_l(z):=E(z-l).
$$
Shifting in this way the function $E$ by a large $l>0$ makes the logarithmic tract is backward invariant and yields $J_{\disF_l}\subset G$. Consequently, only the dynamics of $\disF_l$ in $G$, the domain on which $\disF_l$ is close to the function $\disf_l$, given by the formula 
$$
\disf_l(z):=f(z-l),
$$
is relevant for our purposes. The details of this and the definition of the Julia set in the present setting are given in Section 2.1. 

\bfact[Stallard \cite{Stallard-IV}] Let $p>0$.
All the functions $\disF_l$, $l\in\C$, belong to the Eremenko--Lyubich class $\cB$ and
there exists a constant $C_p>0$ such that for all real $l> C_p$ we have that
$$
\HD(J_{\disF_l})= 1+\frac{1}{1+p} <2\,.
$$
\efact

In the present note we analyze the hyperbolic dimension of these functions. In fact, we first work with the functions $\disf_l$
and then transfer the results to the globally defined entire functions $\disF_l$.

The key point is to employ the thermodynamic formalism
of \cite{MU-integralMean} and, in particular, the Bowen's Formula from this paper that determines hyperbolic dimension.
We will see that $\lim_{l\to \infty}  \hypd (\disF_l) =1$ which clearly implies that $\hypd (\disF_l) < \HD(\J_{\disF_l})$ provided that $l>C_p$ is large 
enough.

\smallskip

\noindent
 {\em Acknowledgement:} 
 We would like to thank the referee for excellent refereeing, comments, and suggestions which improved the final exposition of our results.

\subsection{Notation}

We use standard notation such as 
$\D(z,r)$ for the open disk in $\C$ with center $z\in \C$ and radius $r>0$.
When the center is the origin, we also use the simplified notation
$$
\D_r:=\D(0,r).
$$
The complement of its closure will be denoted by
$$
\D_r^*:=\C\sms \ov\D_r.
$$
Frequently we deal with  half--spaces. Let
$$ 
\cH_s:=\big\{z\in \C: \; \Re {z} >s\big\} \quad , \quad s\geq 0 \, .$$
 When $s=0$, then we also write $\cH$ for $\cH_0$.
 
Many constants, especially those in Fact \ref{fact 2}, depend on the parameter $p$ of the definitions of the functions $E$ and $f$. However, this will be fixed throughout the whole paper and we may ignore it.
 
We say that  
$$
A\preceq B
$$
for non--negative real expressions $A$ and $B$ if and only if there exists a positive constant $C$ independent of variable parameters involved in $A$ and $B$ such that $A\le CB$. We then say that
$A\succeq B$ if and only if $B\preceq A$. Finally, $A\asymp B$ if and only if $A\preceq B$ and $B\preceq A$.

 \section{Singularities, models and approximating entire functions}

\subsection{General definitions}

 Iversen's classification of singularities is explained in length in \cite{KU2}, see also \cite{BE08.1}.
An entire function $g:\C\to\C$ can have only two types of singular values. Firstly, a point $b\in \hat\C$ is a \emph{critical value} of $g$ if and only if $b=g(c)$ for some $c\in \C$ with $g'(c)=0$. 
Secondly, a complex number $b\in \hat \C$ is an \emph{asymptotical value} of $g$ if and only if there exists a continuous function $\ga: [0,+\infty)\to\C$ such that 
$$
\lim_{t\to+\infty}\g(t)=\infty 
\and
\lim_{t\to+\infty}f(\g(t))=b.
$$
In this latter case for every $r>0$ there exists an unbounded connected component $\Om_r$ of $g^{-1}(\D(b,r))$
such that 
$$
\Om_{r'}\subset \Om _r
$$ 
whenever $r'<r$ and 
$$
\bigcap_{r>0} \Om_r =\emptyset.
$$
Such a choice of components is called an asymptotic tract over $b$
and it is called \emph{logarithmic tract} in the case when the map $g:\Om_r\to \D(b,r)\setminus \{b\}$ is a universal covering for some $r>0$. The set of singular values of $f$ is proved to consist of all critical and asymptotic values of $f$. Its intersection with $\C$
will be denoted by $S(g)$.

We consider functions belonging to
 the Eremenko--Lyubich class $\cB$ that consists of all entire functions $g$ for which $S(g)$ is a bounded set. These functions are also called of \emph{bounded type. 
}
If $g\in \cB$, then there exists $r>0$ such that $
S(g)\subset \D_r
$.
Then $g^{-1}(\D_r^*)$ consists of mutually disjoint unbounded Jordan domains $\Omega_r$ with real analytic boundaries such that $g:\Om\to \D^*_r$ is a covering map (see \cite{ EL92}). Thus, an entire function $g$ in class $\cB$
has only logarithmic singularities over infinity. As we already mentioned it, the connected components of $g^{-1}(\D^*_r)$ are called \emph{tracts} or, more precisely, \emph{logarithmic tracts}. Then there exist all holomorphic branches of the logarithm of $g$ restricted to $\Om_r$. Fix one of them and denote it by $\tau$. So, 
\beq\label{Nov20 1}
g|_{\Om_r} =\exp\circ{\tau},
\eeq
where 
$$
\ph=\tau^{-1} :\cH_{\log r} \to \Om_r
$$
is a conformal homeomorphism. In addition, $\ph$ extends continuously to $\infty$ and $\ph(\infty)=\infty$. 

Keeping this notation, if we restrict $g$ to the tracts over infinity then it is now standard, especially since the appearance
of the papers \cite{Bishop-EL-2015, Bishop-S-2016} by Chris Bishop, to call the map
$$
g_{|g^{-1}(\D_r^*)}: g^{-1}(\D_r^*) \to \D_r^*
$$
a model function. We will see that the functions considered in our current paper have only one single tract over infinity. This is the reason why we use the following simplified definition of a model function. This is in the spirit of the definition
in \cite{Rempe-HypDim2},
see \cite{Bishop-EL-2015, Bishop-S-2016} for the general one.

\bdfn\label{tract model}
A model is any  holomorphic map 
$$
g=e^\tau :\Om_r \to \D_r^*,
$$
where 
\ben 
\item $r\in [1,+\infty)$, 

\item $\Om_r$ is a simply connected unbounded domain in $\C$, called a tract, such that $\partial\Om_r$ is a connected subset of $\C$ 

and 
\item 
$\tau : \Tract_r \to \cH_{\log r}$ is a conformal homeomorphism fixing infinity; the latter more precisely meaning that
$$
\tau (z)\to \infty\;\; {\rm as } \;\;z\to \infty. 
$$
\een
\edfn

\medskip

The tract $\Om_r$ may or may not intersect the disk $\D_r$. The later case has important dynamical consequences.

\bdfn\label{defn disj type}
If $f$ is a model or an entire function of bounded type and if there exists $r>0$ such that 
\beq\label{2023-1}
S(f)\subset \D_r \quad  \  \text{and} \quad    \ov {f^{-1}({\D_r^*}) }\sbt \D_r^*,
\eeq
then  $f$ is called of disjoint type. 
\edfn

If $f$ is such a disjoint type model or entire function, then the \emph{Julia set} of $f$ is defined to be
$$
\J_f:=\big\{z\in \D_r^*: f^n(z)\in \D_r^* \; \; \text{for all}\;\; n\geq 1\big\}\,.
$$
For disjoint type entire functions this definition coincides with the usual one, see Proposition~2.2 in \cite{R-ALC}.


\medskip

\subsection{Elementary properties of the functions $E$ and $f$} We now discuss some elementary properties of the functions introduced in the introduction and we examine how they behave with respect to the above definitions. 

To start with, we recall that these functions have been introduced and studied by Stallard and her paper \cite[Section 3]{Stallard-III} lists elementary properties of $E$ and $f$. We recall now some necessary facts from this paper.

Let's denote:
$$
G_{x_0,\ka}:= \left\{z=x+iy\in\C: x > x_0 \text{ and } |y| < \ka \frac{\pi x }{ (1+p)(\log x)^p}\right\}
$$
and abbreviate  $G_{x_0} = G_{x_0,1}$ 
so that the set $G$ of \eqref{eq 3} is $G_3=G_{3,1}$. 
For any integer $n\geq 3$ let $\sg_{n+1}$ be the boundary of the open set $G\setminus \ov G_{n+1}$. The orientation of $\sg_{n+1}$ and of all following boundary curves are always understood  in the clockwise direction. Cauchy's Integral Formula shows that 
$$
 \frac{1}{2i\pi}\int_{\sg_{n+1}} \frac{f(\xi)}{\xi-z} d\xi =0 \quad \text{for every $z\not\in \ov G$.}
$$
Therefore, still for  $z\not\in \ov G$,
$$
E(z)=  \frac{1}{2i\pi}\int_{\partial G} \frac{f(\xi)}{\xi-z} d\xi =  \frac{1}{2i\pi}\int_{\partial G_{n+1}} \frac{f(\xi)}{\xi-z} d\xi\,.
$$
It thus follows that the right hand side integral gives the holomorphic extension of $E$ to the domain 
$\C\setminus \ov G_{n+1}$. 

Consider now an arbitrary point $z\in G\setminus \ov G_{n+1}$. Then, Cauchy's Residue Theorem shows that 
\begin{equation*}
\begin{split}
E(z)= &  \frac{1}{2i\pi}\int_{\partial G_{n+1}} \frac{f(\xi)}{\xi-z} d\xi\\
 =& 
-\frac{1}{2i\pi}\int_{\partial (G_{n}\setminus \ov G_{n+1})} \frac{f(\xi)}{\xi-z} d\xi +  \frac{1}{2i\pi}\int_{\partial G_{n}} \frac{f(\xi)}{\xi-z} d\xi\\
=&
\; f(z)+ \frac{1}{2i\pi}\int_{\partial G_{n}} \frac{f(\xi)}{\xi-z} d\xi.\,
\end{split}
\end{equation*}

Starting with this observation, one can get the following fact which is contained in Lemma 3.1 in \cite{Stallard-III} along with its proof.

\bfact \label{fact 2} Let $\widecheck L , \widehat L$
be the boundary of $G_{D+1, \frac56}, G_{D-1, \frac76}$ respectively.
Then there exist constants $C,D>3$
such that the following hold. 
\ben
\item If $z\not\in G_{D}$ then 
$$
|E(z)|\leq C
$$ 
and

\noindent if $z\in G_{D}$ then 
$$
|E(z)-f(z)|\leq C \  {\rm\text{ as well as }} \ |E'(z)-f'(z)|\leq C.
$$
\item 
$$ 
E(z) =\frac1{2\pi i}\int _{\widecheck L} \frac{f(t)}{t-z}dt \quad \text{ for $ z\not\in G_{D}$}$$
and
$$
E(z)=f(z)+ \frac1{2\pi i}\int _{\widehat L} \frac{f(t)}{t-z}dt \quad \text{ for $ z\in G_{D}$}.
$$
\item If $z\in G_{D,\frac76}\setminus \Int (G_{D,\frac56})$ then 
$$
|f(z)|\leq \exp \Big(-\frac12e^{\frac12(\log \Re z)^{1+p}} \Big).
$$
\een   
\efact

Item (1) from this Fact  \ref{fact 2} shows that $f^{-1}(\D_r^*)\subset G_D$
for  every $r>2C$. Elementary estimates, based on the explicit representation of $f$, show that
$f^{-1}(\D_r^*)$ is a simply connected unbounded domain in $\C$. It turns out that the same is true for the approximating entire function $E$; details can be found in Proposition 2.2 of \cite{Rempe-HypDim2}.
Thus, we have the following.

\bfact\label{21 1}
Let $C$ be given by Fact  \ref{fact 2}. Then, there exists $r_0>4C$ such that 
$$
S(E)\subset \D_{r_0/2} 
$$
and for every $r\geq r_0/2$, both sets $E^{-1}(\D_r^*)$ and $f^{-1}(\D_r^*)$ are simply connected unbounded domains in $\C$ contained in $G_D$. They will be respectively denoted by
$$
\Om_{E,r}:= E^{-1}(\D_r^*) \quad \text{and}\quad \Om_{f,r}:=f^{-1}(\D_r^*)\,.
$$
\efact

\smallskip

\noindent
From now on fix any 
\beq\label{2023-12}
r\geq r_0/2,
\eeq
where $r_0$ comes from Fact \ref{21 1}.
Then the map $f:\Om_{f,r} \to \D_r^*$ is of the form $f(z)=e^{\tau (z)}$ with  $\tau : \Om_{f,r} \to \cH_{\log r}$ 
given by 
$$
\tau(z):=\exp ((\log z)^{1+p}).
$$
We have to know what the inverse conformal homeomorphism
$\ph = \tau^{-1} : \cH_{\log r}\to \Om_{f,r}$ looks like. Indeed, a straightforward calculation gives
\beq\label{ph}
\ph (\xi ) = \exp \left((\log \xi )^{\frac{1}{1+p}}\right),
\eeq
where $\log $ is the principal branch of logarithm again, i.e. determined by the requirement that $\log 1=0$. 

In conclusion,
\beq\label{21 3}
f_{| \Om_{f,r}}=e^\tau:\Om_{f,r} \to \D_r^*
\eeq
 is  a model as defined in Definition \ref{tract model}; fact \ref{fact 2} explains how the entire function $E$ approximates this model. 

\blem\label{Koebe}
There exists a constant $K\geq 1$ such that
$$
\frac1K\leq \frac{|\ph'(\xi+iy)|}{|\ph'(\xi )|}\leq K
$$
 for every $\xi$ with $\Re(\xi) \geq \log r_0$ and every $0\leq y\leq 2\pi$.
\elem

\bpf
The statement follows from Koebe's Distortion Theorem since the conformal map $\ph = \cH_{\log r_0}\to \Om_{f,r_0}$ is in fact defined on the half space $\cH_{\log(r_0/2)}$.
\epf

\subsection{Disjoint Type Versions of $E_l$ and $f_l$}
Given any $l\in\C$, the functions $\disf_l =f\circ T_l$ and $ \disF_l=E\circ T_l$, where $T_l$ is the translation $z\mapsto z-l$, have been defined in the introduction. We have that $ \disF_l\in \cB$ since it is known, see \cite{Stallard-IV}, that $E\in \cB$.
  
Obviously,
\beq\lab{120230411}
\Om_{\disf_l,r}:=\disf_l ^{-1} (\D_r^*) = f^{-1} (\D_r^*) +l = \Om_{f,r}+l,
\eeq
and also 
\beq\lab{220230411}
\Om_{\disF_l,r}=\disF_l^{-1} (\D_r^*) = E^{-1} (\D_r^*) +l.
\eeq
By Fact \ref{21 1}, for all $r\ge r_0/2$ and $l\in[0,+\infty)$, all these tracts are contained in respective sets $G_D+l$. So, setting 
\beq\label{def l0} l_r:= \max\{0, r-D\}\, ,\eeq
we have that
\beq\label{2023-5}
\Om_{\disf_l,r}  \;, \;  \Om_{\disF_l,r} \subset \D_{r}^* 
\eeq
for all $r\ge r_0/2$ and  all $l\geq l_r$.
Consequently, all the functions $\disf_l , \disF_l$, $l\geq l_r$, are of disjoint type and for their Julia sets we have that
\beq\label{2023-6}
\J_{\disf_l} , \J_{\disF_l}   \subset \D_{r}^*
\eeq
for all $r\ge r_0/2$ and  all $l\geq l_r$.

Recall that for the model $f$ we have the expression
\eqref{21 3}. The analogous expression for $\disf_l$ is 
\beq\label{2023-9}
{\disf_l}_{|\Om_{\disf_l ,r}} =e^{\tau_l} : \Om_{\disf_l,r} \to \D_r^*
\eeq
where $\tau_l (z) =\tau (z-l)$ so that the inverse of $\tau_l$ is 
\beq\label{2023-10}
\ph_l =\ph +l : \cH_{\log r} \to \Om_{\disf_l, r}
\eeq
where $\ph$ is still the conformal map defined by \eqref{ph}

\section{Thermodynamical formalism}
Our ultimate goal is to determine the hyperbolic dimension of the functions $\disF_l$ which, under certain conditions, can be done by employing the methods of thermodynamic formalism.
The hyperbolic dimension is then given by the zero of the topological pressure, the fact that goes back to Bowen \cite{Bow79}. In the present context, namely for disjoint type
models and entire functions of bounded type, such a theory has been developed in  \cite{MU-integralMean}.

Let $\cC_b (\D_r^*)$ be the vector space of all complex--valued bounded continuous functions defined on $ \D_r^*$. Endowed with  the supremum norm, it becomes a Banach space.

Let $g:=\disf_l$ or $g:=\disF_l$. Given $t > 0$, the transfer operator for the map $g$ and for the parameter $t$, acting on a function $h\in \cC_b (\D_r^*)$, is defined by the formula 
\beq\label{transfer operator}
\pf_{g,t} h (w) := \sum_{g(z)=w} |g'(z)|^{-t}_1 h(z) \  \  \text{for every}\  w\in \D_r^*,
\eeq
where 
$$
|g '(z)|_1:= \frac{|g'(z)|}{|g(z)|}|z| 
$$ 
is the logarithmic derivative of $g$ evaluated at the point $z$.
 
We are to find out for which parameters $t>0$ the following two crucial properties hold:
  \beq\label{crucial} \|\pf_{g,t} \1 \|_\infty<+\infty \quad \text{and} \quad \lim_{w\to\infty} \pf_{g,t}\1 (w) =0.
  \eeq
 

Indeed, since our map $g$ is of disjoint type, once \eqref{crucial} is verified  then, following \cite[Section 8]{MU-integralMean}, we deduce that the whole thermodynamic formalism, along with all its applications obtained in \cite{MU-integralMean}, holds. Especially Bowen's Formula does. This formula involves topological pressure which for the disjoint type map $g$ is given at a parameter $t\in (0,+\infty)$ by the formula
 \beq\label{topo pressure}
 \P(g,t)=\lim_{n\to\infty} \frac1n\log  \pf_{g,t}^n\1(w),
 \eeq
 where $w\in \D_r^*$ is any arbitrarily chosen point. The limit  exists and is independent of $w$ because of Theorem~8.1 in \cite{MU-integralMean} which ultimately goes back to Lemma~5.8 and Corollary~5.18 in \cite{MU-Memoir}.
 
 \smallskip
 
  \subsection{Estimates for the Transfer operators of the Model Functions $\disf_l$}
 
 \bprop\label{9.2-1}
 Let $\pf_{\disf_l,t}$ be the transfer operator of $\disf_l$, $l\geq 0$, with a parameter $t>0$. Fix $r\ge r_0$. Let $w_0\in \D_r^*$. Then
 $$
 \pf_{\disf_l,t} \1 (w_0) <\infty \quad \text{if and only if} \quad t>1\,.
 $$
Moreover, if $t>1$ then \eqref{crucial} holds for $g=\disf_l$.
 \eprop
 
\bpf
Having $w_0\in \D_r^*$ and $t>0$, let us start exactly as in  the proof of Theorem~4.1 in \cite{MU-integralMean}.
 If $z_l\in \disf_l^{-1}(w_0)$ then, using \eqref{2023-9}, the logarithmic derivative can be expressed as follows:
$$
|\disf_l '(z_l)|_1 = |\tau_l'(z_l) z_l| = \frac{|\ph_l (\xi)|}{|\ph'_l (\xi)|} = |(\log \ph_l )'(\xi)|^{-1}
$$
where $\xi =\tau_l(z_l)$ and where $\ph_l=\ph+l$ is the map of
\eqref{2023-10}. Notice that $\xi=u+iv$ does not depend on $l$, where $u=\log|w_0|$. From this, together with Lemma~\ref{Koebe}, we get that
$$
\pf_{\disf_l,t}\1 (w_0)=\sum_{\exp(\xi )=w_0}  |(\log \ph_l )'(\xi)|^{t} \asymp \int _\R  |(\log \ph_l )'(\log|w_0|+iv)|^{t} dv\,.
$$
Now, since $\ph_l=\ph+l$ and since we have the explicit expression \eqref{ph} for  $\ph $, we can calculate as follows:
$$
|(\log \ph_l )'(\xi)|= \left| \frac{\ph(\xi )}{\ph(\xi)+l}\right| \frac1{1+p} \frac1{|\xi| |\log \xi |^{\frac{p}{1+p}}}\asymp  \left| \frac{\ph(\xi )}{\ph(\xi)+l}\right| \frac1{|\xi|(\log |\xi |)^{\frac{p}{1+p}}}.
$$
since $\arg (\xi)\in (-\pi/2, \pi/2)$. Therefore,
\beq\label{estimate transfer}
\pf_{\disf_l,t}\1 (w_0) \asymp \int _\R  \left| \frac{\ph(\xi )}{\ph(\xi)+l}\right|^t \frac1{|\xi|^t(\log |\xi |)^{\frac{tp}{1+p}}}dv.
\eeq
Since $\lim_{|v|\to+\infty}\ph (\log|w_0|+iv)=\infty$, we have that
$$
\frac23 \leq \left| \frac{\ph(\xi )}{\ph(\xi)+l}\right|\leq 2
$$
whenever $|v|=|\Im (\xi)|$ is sufficiently large.  Thus we get from \eqref{estimate transfer}
that $\pf_{\disf_l,t}\1 (w_0)$ is finite if and only if $t>1$.

\sp The uniform bound of $\| \pf_{\disf_l,t}\|_\infty <\infty$ also follows from \eqref{estimate transfer}. Indeed, let $w=e^\xi \in \D_r^*$. Then $z=\ph(\xi)\in G_D$, whence $x=\Re (z)>0$. Thus,
\beq\label{ineg ridiculus}
\left| \frac{\ph(\xi )}{\ph(\xi)+l}\right|^2=\frac{x^2+y^2}{(x+l)^2+y^2}\leq 1.
\eeq
It follows from this that 
$$
\pf_{\disf_l,t}\1 (w) \preceq \int _\R \frac1{|\xi|^t(\log |\xi |)^{\frac{tp}{1+p}}}dv 
=\frac12
\int _\R \frac1{(u^2+v^2)^{\frac{t}{2}}(\log (u^2+v^2))^{\frac{tp}{1+p}}}dv .
$$
Since for every for $w\in \D_r^*$ we have $u\geq u_r=\log r$ it follows that
\beq\lab{520230404}
\sup_{w\in \D_r^*}\pf_{\disf_l,t}\1 (w) 
\preceq C:=\int_\R \frac1{(u_r^2+v^2)^{\frac{t}{2}}(\log (u_r^2+v^2))^{\frac{tp}{1+p}}}dv<+\infty .
\eeq
\smallskip

Finally, if $t>1$ then $\d =(t-1)/2>0$, whence 
\beq\label{uniformly to zero}
\pf_{\disf_l,t}\1 (w) \preceq \frac1{u^\d}\int _\R \frac1{|u_r+iv|^{1+\d} }dv \preceq \frac1{(\log |w|)^\d}.
\eeq
This shows that $\lim_{w \to\infty} \pf_{\disf_l,t} \1 (w)=0$.
\epf
 
 The next result gives an estimate for the topological pressure. More precisely, it shows that for a given $t>1$ the pressure $\P(\disf_l ,t)<0$ for all sufficiently large values of $l$.
  
\bprop\label{prop pressure model} Let $t>1$. Fix $r\ge r_0$. Then,
 for every $\e>0$ there exists $l_{\e,r,t} \geq l_r$ such that 
 $$
 \pf _{\disf_l ,t} \1 (w) \leq \e
 \quad \text{for every $l\geq l_{\e,r,t}$ and every $w\in \D_r^*$}
  \,.
 $$
 \eprop
 
 \bpf
 Let $t>1$ and $\e>0$. We are in the same situation as in the proof of Proposition \ref{9.2-1}.
 The first benefit we take out of this proof is that the convergence  $\lim_{w\to\infty}\pf_{\disf_l ,t}\1 (w) =0$
is uniform in $l\geq 0$; see \eqref{uniformly to zero}. Therefore,
there exists $r_\e \ge r$ such that 
 $$
 \pf_{\disf_l ,t}\1 (w) <\e \quad \text{whenever $|w|\geq r_\e$ and $l\geq 0$.}
 $$
Moreover, this proof shows that the integral
 $$
 \int _\R   \frac1{|\xi|^t(\log |\xi |)^{\frac{tp}{1+p}}}dv  \, \quad \xi =u+iv\, ,
 $$
 converges uniformly for $u\ge u_r=\log r$. Therefore, there exists $V=V_{\e,t}$ such that
$$
 \int _{|v|\geq V}  \frac1{|\xi|^t(\log |\xi |)^{tp/1+p}}dv  \leq \frac{\e}{2} \quad \text{for every $u\ge u_r$.}
$$
So, by invoking now \eqref{estimate transfer} and \eqref{ineg ridiculus}, we conclude that it remains to estimate the integral
$$
 \int _{|v| < V}    \left| \frac{\ph(\xi )}{\ph(\xi)+l}\right|^t \frac1{|\xi|^t(\log |\xi |)^{\frac{tp}{1+p}}}dv
$$
from above by $\e/2$ for all $l\ge 0$ large enough and all $w\in \D_r^*\setminus \D_{r_\e}^*$. Here we used again the notation $\xi =u+iv$, $u=\log |w|$. Notice that all points $\xi$ that appear in this integral belong to the compact set
 $$
 K =\{\xi=u+iv:\log r \leq u \leq \log r_\e \text{ and } |v|\leq V\}.
 $$
Since $M:=\sup_{\xi\in K} \{|\ph (\xi )|\}<+\infty$, we have that
$$
\left|\frac{\ph(\xi )}{\ph(\xi)+l}\right|\leq \frac{M}{l-M} \quad , 
$$
for every $l>M$ and all $\xi \in K$. Thus,
$$
 \int _{|v| < V}\left| \frac{\ph(\xi )}{\ph(\xi)+l}\right|^t \frac1{|\xi|^t(\log |\xi |)^{\frac{tp}{1+p}}}dv\leq C  \frac{M}{l-M}
 \le \frac{\e}{2},
$$
where $C\in(0,+\infty)$ is the constant coming from \eqref{520230404} and the last inequality was written assuming that $l$ is large enough.
\epf
 
\smallskip
 
\subsection{Behavior of the Transfer Operators for Entire functions $\disE_l$}
We now have sufficiently strong estimates for the transfer operators of the models $f_l$. Since ultimately we are after the entire functions $\disF_l$,
we have to carry over these estimates to the transfer operators of these functions $\disF_l$. Since the entire functions approximate the models, i.e. since we have Fact \ref{fact 2}, we are in a similar situation as in \cite{MyZdunik2021} where also the operators of some models and approximating entire functions have been compared. Following the approach of that paper we will prove the following.
 
\bprop\label{2023-11} 
There exist constants $\cK\in[1,+\infty)$ and $r_1\ge r_0$ such that for every $t>1$, all $l\in\C$, and all $r\geq r_1$, we have that
$$
\frac1{\cK^t}\leq \frac{\pf_{\disF_l,t}\1 (w)}{\pf_{\disf_l,t}\1 (w)}\leq \cK^t\quad \text{for all} \quad w\in \D_r ^*.
$$
\eprop

In our proof of Proposition~\ref{2023-11}  we adapt here the approach of \cite{MyZdunik2021}, particularly Section 7 of that paper.
We will show that \cite[Lemma 7.3]{MyZdunik2021}
holds in the present setting if $r\geq r_0$ is large enough. This will suffice. We first shall prove the following.

\bfact\label{2.3} For all sufficiently large $r\geq r_0$, say $r\ge r_1\ge r_0$, we have that
$$
\frac12 \leq \frac{|\disF_l(z)|}{|\disf_l(z)|}\leq 2 \quad \text{and}\quad \frac12 \leq \frac{|\disF'_l(z)|}{|\disf'_l(z)|}\leq 2
$$
for all $l\in\C$ and all $z\in \Om_{\disf_l,r}$.
\efact

\bpf
The first inequality is a direct consequence of item (1) in Fact \ref{fact 2} combined with the inequality $r\geq r_0> 4C$ established in Fact~\ref{21 1}.

In order to proof the second inequality we also start with item (1) in Fact \ref{fact 2}. It gives
$$
\left|\frac{|E'(z)|}{|f'(z)|}-1\right|\leq \frac{C}{|f'(z)|} \quad \text{for all $z\in \Om_{f,r}\subset G_D$.}
$$
This time we have to estimate $|f'(z)|$ and to show that there exists some $r \geq r_0$ such that 
\beq\label{2.5}
\frac{C}{|f'(z)|}\leq \frac 12 \quad \text{for all $z\in \Om_{f,r}$.}
\eeq
Remember that $f(z)=e^{\tau(z)}=e^{\ph^{-1}(z)}$ for every $z\in \Om_{f,r}$. Thus,
$$
f'(z)=\frac{f(z)}{\ph'(\xi )}\quad\text{where $\xi = \ph^{-1}(z)\in \cH_{\log r}$.}
$$
Obviously  $|f(z)|>r$ but what about $|\ph'(\xi )|$? From the formula \eqref{ph} we get
$$
\ph'(\xi) = \frac1{1+p}\exp ((\log \xi)^{\frac{1}{1+p}})\frac1{(\log \xi)^{\frac{p}{1+p}}\,\xi} .
$$
If $v:=\log \xi$ then
\beq\label{2.4}
|\ph'(\xi)| 
\leq \lt|\frac{\exp (v^{\frac{1}{1+p}})}{v^{\frac{p}{1+p}}e^v}\rt| =\frac{\exp\lt(\Re\(v^{\frac{1}{1+p}}-v\)\rt)}{|v|^{\frac{p}{1+p}}}.
\eeq
Since $\xi\in \cH_{\log r}$, $\Re v > \log \log r $, and $|\Im v|<\pi/2$, so if we write $v=se^{i\al}$, then
$$
s> \log\log r \quad \text{and}\quad |\al| < \frac{\pi/2}{\log\log r} .
$$
Thus,
$$
\Re\(v^{\frac{1}{1+p}}-v\)
=-s\lt(\cos \al -s^{-\frac{p}{1+p}}\cos\Big(\frac{\al}{1+p}\Big)\rt)
\leq -\frac s2
\leq -\frac{\log\log r}{2}
$$
 provided $r$ is sufficiently large. In this case we get from \eqref{2.4} that
$$
|\ph'(\xi)| \leq \frac{1}{\sqrt{\log r }(\log \log r)^{\frac{p}{1+p}}}.
$$
This shows that \eqref{2.5} holds for all $r\ge r_0$ sufficiently large. Thus, \ref{2.3} holds for $f$ and $E$, i.e. if $l=0$. It then holds for all $l\in\C$ because of \eqref{120230411}.

\epf

Having established Fact \ref{2.3}, the proof of Proposition 7.4 in \cite{MyZdunik2021} applies word by word and shows that
the required inequality in
Proposition \ref{2023-11} holds.

\medskip

\section{Proof of Theorem \ref{thm main}}
As it was explained in the Introduction, it suffices to show that
\beq\lab{120230412}
\lim_{l\to \infty}  \hypd (\disF_l) =1.
\eeq
In order to do this fix $t>1$. Fix also any $r\ge r_1$, for example $r=r_1$. By virtue of Proposition~\ref{prop pressure model}, we have that
$$
\pf_{f_l,t}\1(w)\le \cK^{-t}
$$
for all $l\ge l_{\cK^{-t},r,t}$ and all $w\in\D_r^*$. So, by Proposition~\ref{2023-11},
$$
\pf_{E_l,t}\1(w)\le 1
$$
for all $l\ge l_{\cK^{-t},r,t}$ and all $w\in\D_r^*$. In conjunction with \eqref{topo pressure}, this gives that
$$
\P(E_l,t)\le 0
$$
for all $l\ge l_{\cK^{-t},r,t}$. So, if $X\sbt \J_{E_l}$ is an arbitrary hyperbolic set for $E_l$, then 
$$
\P(E_l|_X,t)\le 0.
$$
The supremum over all hyperbolic sets of the left hand side of this inequality is the hyperbolic pressure
$P_{hyp}(E_l , t)$ of $E_l$ evaluated at $t$. So, we have that
\beq \label{hyp pressure neg}
\P_{hyp}(E_l , t)=\sup \{ \P(E_l|_X,t):\text{$X$ is a hyperbolic set for $E_l$}\}\leq 0.
\eeq
Now, we want to use the Bowen's Formula of \cite{BKZ-2009}. Theorem B of this paper applies to the functions $E_l$
and states that the hyperbolic dimension of the set $E_l$ is equal to
$$
\hypd (\disF_l)=\inf \{s>0: \P_{hyp}(E_l , s) \leq 0 \}.
$$
Combined with \eqref{hyp pressure neg}, we thus get that 
$$
\hypd (\disF_l)\le t
$$
for all $l\ge l_{\cK^{-t},r,t}$. So, the formula \eqref{120230412} is established and the proof of Theorem \ref{thm main} is complete. 

\medskip


\bibliographystyle{plain}



\end{document}